\documentclass[11pt]{article}
\usepackage{amsmath}
\usepackage{amssymb}
\usepackage{amsfonts}
\usepackage{eucal}
\usepackage{graphicx,url}
\usepackage{amssymb,amsmath, amsthm}
\usepackage{epstopdf}
\input pictex.sty
\newtheorem{theorem}{theorem}[section]
\newtheorem{proposition}[theorem]{Proposition}

\newtheorem{example}{Example}
\def\bincoeff#1#2{{#1\choose #2}}
\begin{document}
\title{Fractional discrete processes: compound and mixed Poisson representations}
\author{Luisa Beghin\thanks{Address: Dipartimento di
Scienze Statistiche, Sapienza Universit\`a di Roma, Piazzale Aldo
Moro 5, I-00185 Roma, Italy. e-mail:
\texttt{luisa.beghin@uniroma1.it}}\and Claudio
Macci\thanks{Dipartimento di Matematica, Universit\`a di Roma Tor
Vergata, Via della Ricerca Scientifica, I-00133 Roma, Italia.
E-mail: \texttt{macci@mat.uniroma2.it}}}
\date{}
\maketitle
\begin{abstract}
\noindent We consider two fractional versions of a family of
nonnegative integer valued processes. We prove that their
probability mass functions solve fractional Kolmogorov forward
equations, and we show the overdispersion of these processes. As
particular examples in this family, we can define fractional
versions of some processes in the literature as the Polya-Aeppli,
the Poisson Inverse Gaussian and the Negative Binomial. We also
define and study some more general fractional versions with two
fractional parameters.
\end{abstract}
\maketitle

\section{Introduction}\label{sec:introduction}
In this paper we consider a large family of nonnegative integer
valued processes $\{M(t):t\geq 0\}$ defined by
\begin{equation}\label{eq:cs0}
M(t):=\sum_{k=1}^{N_\lambda(t)}X_k,
\end{equation}
where $\{X_n:n\geq 1\}$ is a sequence of i.i.d. positive and
integer valued random variables, independent of a (non-fractional)
Poisson process $\{N_\lambda(t):t\geq 0\}$ with intensity
$\lambda$. Throughout this paper we also deal with a mixed Poisson
representation of the process in \eqref{eq:cs0}; more precisely we
mean $\{N_1(S(t)):t\geq 0\}$, where $\{S(t):t\geq 0\}$ is a
subordinator, independent of a (non-fractional) Poisson process
$\{N_1(t):t\geq 0\}$ with intensity $1$, such that $N_1(S(t))$ is
distributed as $M(t)$ for each fixed $t>0$.

Our aim is to present and analyze two fractional versions of the
process $\{M(t):t\geq 0\}$ (see \textbf{(FV1)} and \textbf{(FV2)}
at the beginning of Section \ref{sec:representations}). This
allows to define fractional generalizations of some processes in
the literature, which include the Polya-Aeppli, the Poisson
Inverse Gaussian and the Negative Binomial. These processes are
commonly used when the empirical count data exhibit
overdispersion, i.e. when the sample variance is larger than the
sample mean; moreover it is known that mixed Poisson processes
provide simple counting models with overdispersion. We shall see
that the fractional versions have the same feature.

The two fractional versions are obtained by considering
independent random time-changes of $\{M(t):t\geq 0\}$ in terms of
a stable subordinator $\mathcal{A}^\alpha$ or its inverse
$\mathcal{L}^\alpha$ (for $\alpha\in(0,1)$). This approach is
inspired by the recent increasing interest on random time-changed
and subordinated processes; see, among the others, \cite{Sato},
\cite{MeerschaertScheffler}, \cite{HeydeLeonenko}, \cite{LindeShi}
and \cite{KumarNaneVellaisamy}. These processes are widely studied
and applied, mainly because their finite dimensional distributions
display power law decay and thus heavy tails, either in the case
of renewals and birth processes (see \cite{MainardiGorenfloVivoli}
and \cite{OrsingherPolito2010}) or in that of diffusions (see
\cite{AlloubaZheng}, \cite{SaichevZaslavsky} and
\cite{OrsingherBeghin}). We recall that fractional versions of the
compound Poisson process (as the ones in this paper) have been
studied in \cite{Scalas} and \cite{BeghinMacci}, in the case of
continuous distributed summands.

In this paper we provide a natural fractional extension of some
results in the literature. In particular Proposition
\ref{prop:Kolmogorov-equations} below shows that the probability
mass functions of the two fractional versions of the process in
\eqref{eq:cs0} solve suitable Kolmogorov equations, where the
classical derivatives are replaced by the fractional derivative in
the Caputo sense (for the first version), and by the right sided
fractional Riemann-Liouville derivatives on $\mathbb{R}_+$ (for
the second version).

We conclude with the outline of the paper. We start with some
preliminaries in Section \ref{sec:preliminaries}. In Section
\ref{sec:representations} we illustrate some properties of the
fractional versions of the process in \eqref{eq:cs0}. Some
examples are presented in Section \ref{sec:examples}. Finally
Section \ref{sec:generalized-examples} is devoted to more general
fractional versions of some processes in Section
\ref{sec:examples} (the Polya Aeppli and the Poisson Inverse
Gaussian) with two fractional parameters.

\section{Preliminaries}\label{sec:preliminaries}
We start with some preliminaries on fractional calculus, i.e. we
give the definitions of two fractional derivatives for real
functions defined on $[0,\infty)$.
\begin{enumerate}
\item If $\nu\in(0,1)$, the Caputo derivative of order $\nu$
(see e.g. (2.4.17) in \cite{KilbasSrivastavaTrujillo} with $a=0$)
is defined by
$${}_CD_{0+,t}^\nu f(t):=\frac{1}{\Gamma(1-\nu)}\int_0^t(t-s)^{-\nu}\frac{d}{ds}f(s)ds\ (\mbox{for all}\ t\geq 0).$$
\textbf{Remark.} If we set $\nu=1$, ${}_CD_{0+,t}^\nu$ coincides
with the classical derivative $\frac{d}{dt}$; see e.g. Theorem
2.1(b) in \cite{KilbasSrivastavaTrujillo}.
\item If $\nu\in(1,\infty)$, the right sided fractional
Riemann-Liouville derivative on $\mathbb{R}_+$ of order $\nu$ (see
e.g. (2.2.4) in \cite{KilbasSrivastavaTrujillo}) is defined by
$${}_{RL}D_{-,t}^\nu f(t):=\frac{1}{\Gamma(m-\nu)}\left(-\frac{d}{dt}\right)^m\int_t^\infty\frac{f(s)}{(s-t)^{\nu-m+1}}ds\ (\mbox{for all}\ t\geq 0),
\ \mathrm{where}\ m:=\lfloor\nu\rfloor+1.$$ \textbf{Remark.} If we
set $\nu=1$, ${}_{RL}D_{-,t}^\nu$ coincides with $-\frac{d}{dt}$
(i.e. the opposite of the classical derivative); see e.g. (2.2.5)
in \cite{KilbasSrivastavaTrujillo}.
\end{enumerate}
Hereafter, for simplicity, we always write
$$\frac{d_C^\nu}{dt^\nu}:={}_CD_{0+,t}^\nu\ (\mathrm{for}\
\nu\in(0,1))\ \mathrm{and}\
\frac{d_{RL}^\nu}{dt^\nu}:={}_{RL}D_{-,t}^\nu\ (\mathrm{for}\
\nu\in(1,\infty)).$$

Throughout the paper we often use the symbol $Z(\cdot)$ to mean a
process $\{Z(t):t\geq 0\}$. In view of what follows we recall some
preliminaries on the stable subordinator
$\mathcal{A}^\alpha(\cdot)$ of order $\alpha\in(0,1)$, and its
inverse $\mathcal{L}^\alpha(\cdot)$. More precisely let
$\mathcal{A}^\alpha(\cdot)$ be the L\'{e}vy process (starting at
the origin) such that, for each fixed $t>0$, we have
\begin{equation}\label{eq:mgf-stable-distribution}
\mathbb{E}[e^{\theta\mathcal{A}^\alpha(t)}]=\left\{\begin{array}{ll}
\exp(-(-\theta)^\alpha t)&\ \mathrm{if}\ \theta\leq 0\\
\infty&\ \mathrm{if}\ \theta>0;
\end{array}\right.
\end{equation}
thus, by referring to \cite{SamorodnitskyTaqqu}, the random
variable $\mathcal{A}^\alpha(t)$ has stable distribution of index
$\alpha$ and parameters $\mu=0$, $\theta=1$ and
$\sigma=(t\cos(\frac{\pi\alpha}{2}))^{1/\alpha}$. Furthermore
$\mathcal{L}^\alpha(\cdot)$ is defined by
$$\mathcal{L}^\alpha(t):=\inf\{z:\mathcal{A}^\alpha(z)>t\};$$
then we have $P(\mathcal{L}^\alpha(t)\leq
z)=P(\mathcal{A}^\alpha(z)\geq t)$ for all $z,t>0$ and
\begin{equation}\label{eq:Laplace-basic}
\mathbb{E}[e^{-\theta\mathcal{L}^\alpha(t)}]=E_{\alpha,1}(-\theta
t^\alpha)\ \mbox{for all}\ \theta\geq 0.
\end{equation}

Now we recall the definition of the two following fractional
Poisson processes (see e.g. \cite{MeerschaertNaneVellaisamy} for
$\nu\in(0,1)$ and \cite{OrsingherPolito2012} for the case
$\nu\in(1,\infty)$). In both cases we consider a random
time-change of a non-fractional Poisson process
$N_\lambda(\cdot)$.
\begin{enumerate}
\item For $\nu\in(0,1)$,
$N_\lambda^\nu(\cdot):=N_\lambda(\mathcal{L}^\nu(\cdot))$ where
$N_\lambda(\cdot)$ and $\mathcal{L}^\nu(\cdot)$ are independent.
\item For $\nu\in(1,\infty)$,
$\hat{N}_\lambda^\nu(\cdot):=N_\lambda(\mathcal{A}^{1/\nu}(\cdot))$
where $N_\lambda(\cdot)$ and $\mathcal{A}^{1/\nu}(\cdot)$ are
independent.
\end{enumerate}
We can give some further details and we need to introduce some
notation. For all integer $r\geq 0$ and for all
$\gamma\in\mathbb{R}$, the rising factorial, also called
Pochhammer symbol, is defined by
$$(\gamma)^{(r)}:=\left\{\begin{array}{ll}
\gamma(\gamma+1)\cdots (\gamma+r-1)&\ \mathrm{if}\ r\geq 1\\
1&\ \mathrm{if}\ r=0,
\end{array}\right.$$
the falling factorial is defined by
$$(\gamma)_r:=\left\{\begin{array}{ll}
\gamma(\gamma-1)\cdots (\gamma-(r-1))&\ \mathrm{if}\ r\geq 1\\
1&\ \mathrm{if}\ r=0,
\end{array}\right.$$
and we also consider the notation
$$\bincoeff{\gamma}{r}:=\left\{\begin{array}{ll}
\frac{(\gamma)_r}{r!}=\frac{\gamma(\gamma-1)\cdots(\gamma-r+1)}{r!}&\ \mathrm{if}\ r\geq 1\\
1&\ \mathrm{if}\ r=0.
\end{array}\right.$$
\begin{enumerate}
\item If $\nu\in(0,1)$, it is known (see
\cite{MeerschaertNaneVellaisamy}) that $N_\lambda^\nu(t)$ is
distributed as $\sum_{n\geq 1}1_{\{T_1+\cdots+T_n\leq t\}}$, where
$\{T_n:n\geq 1\}$ are i.i.d. random variables with Mittag-Leffler
distribution, i.e. with continuous density
$$f(t)=\lambda t^{\nu-1}E_{\nu,\nu}(-\lambda t^\nu)1_{(0,\infty)}(t)$$
(see e.g. \cite{BeghinOrsingher2009} and
\cite{MainardiGorenfloScalas}), where
$E_{\alpha,\beta}(x):=\sum_{r\geq 0}\frac{x^r}{\Gamma(\alpha
r+\beta)}$ is the Mittag-Leffler function (see e.g.
\cite{Podlubny}, page 17). Then, if we consider the generalized
Mittag-Leffler function $E_{\alpha,\beta}^\gamma(x):=\sum_{r\geq
0}\frac{(\gamma)^{(r)}x^r}{r!\Gamma(\alpha r+\beta)}$, we have
(see formula (2.5) in \cite{BeghinOrsingher2010})
\begin{equation}\label{eq:fpp-density-nu-leq-1}
P(N_\lambda^\nu(t)=k)=(\lambda t^\nu)^k E_{\nu,\nu
k+1}^{k+1}(-\lambda t^\nu)\ \mbox{for all integer}\ k\geq 0.
\end{equation}
\item If $\nu\in(1,\infty)$, we have
\begin{equation}\label{eq:fpp-density-nu-geq-1}
P(\hat{N}_\lambda^\nu(t)=k)=\frac{(-1)^k}{k!}\sum_{r=0}^\infty\frac{(-\lambda^{1/\nu}t)^r}{r!}(r/\nu)_k\
\mbox{for all integer}\ k\geq 0.
\end{equation}
We recall that this process is presented in
\cite{OrsingherPolito2012} by referring to the fractional
difference operator $(1-B)^\alpha$, for $\alpha\in(0,1]$, and it
is called space-fractional Poisson process. More precisely the
probability mass function in \eqref{eq:fpp-density-nu-geq-1}, with
$\nu=1/\alpha$, satisfies the following equations:
\begin{equation}\label{eq:OP}
\left\{\begin{array}{ll}
\frac{d}{dt}P(\hat{N}_\lambda^{1/\alpha}(t)=k)=-\lambda^\alpha (1-B)^\alpha P(\hat{N}_\lambda^{1/\alpha}(t)=k)\ \mbox{for all integer}\ k\geq 0\\
P(\hat{N}_\lambda^{1/\alpha}(0)=0)=1,\
P(\hat{N}_\lambda^{1/\alpha}(0)=k)=0\ \mbox{for all integer}\
k\geq 1,
\end{array}\right.
\end{equation}
where $B$ is the so called \emph{backward shift operator} defined
by $Bf(k)=f(k-1)$ and $B^{r-1}Bf(k)=f(k-r)$; thus in particular we
have
$$(1-B)^\alpha f(k)=\sum_{j=0}^\infty(-1)^j\bincoeff{\alpha}{j}f(k-j).$$
We remark that formula \eqref{eq:fpp-density-nu-geq-1} here
coincides with formula (1.2) in \cite{OrsingherPolito2012}, though
it is written in a slightly different way.
\end{enumerate}

\section{The fractional versions: compound and mixed representations}\label{sec:representations}
We start with the \emph{compound representation} of the fractional
versions of the process $M(\cdot)$ in \eqref{eq:cs0}, and we refer
to the processes defined in Section \ref{sec:preliminaries}:
\begin{quotation}
\noindent \textbf{(FV1)}:
$M^\nu(t):=M(\mathcal{L}^\nu(t))=\sum_{k=1}^{N_\lambda^\nu(t)}X_k$,
for $\nu\in(0,1)$;\\
\ \\
\textbf{(FV2)}:
$\hat{M}^\nu(t):=M(\mathcal{A}^{1/\nu}(t))=\sum_{k=1}^{\hat{N}_\lambda^\nu(t)}X_k$,
for $\nu\in(1,\infty)$.
\end{quotation}
In view of what follows we consider the following notation for the
probability mass function and the probability generating function
of the random variables $\{X_n:n\geq 1\}$:
$$q_k:=P(X_1=k)\ (\mbox{for all integer}\ k\geq 1);\quad g_q(u):=\sum_{k=1}^\infty u^kq_k.$$
We also consider the notation
\begin{equation}\label{eq:convolution-densities}
q_k^{*n}:=P(X_1+\cdots+X_n=k),
\end{equation}
and therefore we have $q_k^{*n}=0$ for all integer $k<n$.

Now we introduce the two following probability mass functions:
\begin{enumerate}
\item for $\nu\in(0,1)$, we have
\begin{equation}\label{eq:compound-densities-nuleq1}
p_k^\nu(t):=P(M_\lambda^\nu(t)=k)=\left\{\begin{array}{ll}
P(N_\lambda^\nu(t)=0)&\ \mathrm{if}\ k=0\\
\sum_{n=1}^kq_k^{*n}P(N_\lambda^\nu(t)=n)&\ \mathrm{if}\ k\geq 1;
\end{array}\right.
\end{equation}
\item for $\nu\in(1,\infty)$, we have
\begin{equation}\label{eq:compound-densities-nugeq1}
\hat{p}_k^\nu(t):=P(\hat{M}_\lambda^\nu(t)=k)=\left\{\begin{array}{ll}
P(\hat{N}_\lambda^\nu(t)=0)&\ \mathrm{if}\ k=0\\
\sum_{n=1}^kq_k^{*n}P(\hat{N}_\lambda^\nu(t)=n)&\ \mathrm{if}\
k\geq 1.
\end{array}\right.
\end{equation}
\end{enumerate}
The aim of this section is to prove that these probability mass
functions satisfy suitable versions of Kolmogorov equations with
fractional derivatives.

\begin{proposition}\label{prop:Kolmogorov-equations}
If $\nu\in(0,1)$, the probability mass function in
\eqref{eq:compound-densities-nuleq1} satisfies the equations
\begin{equation}\label{eq:system1-nuleq1}
\left\{\begin{array}{ll}
\frac{d_C^\nu}{dt^\nu}p_0^\nu(t)=-\lambda p_0^\nu(t)\\
\frac{d_C^\nu}{dt^\nu}p_k^\nu(t)=-\lambda
p_k^\nu(t)+\lambda\sum_{i=1}^kq_ip_{k-i}^\nu(t)\ \mbox{for all
integer}\ k\geq 1,
\end{array}\right.
\end{equation}
with the initial conditions $p_0^\nu(0)=1$ and $p_k^\nu(0)=0$ for
all integer $k\geq 1$.\\
If $\nu\in(1,\infty)$, the probability mass function in
\eqref{eq:compound-densities-nugeq1} satisfies the equations
\begin{equation}\label{eq:system1-nugeq1}
\left\{\begin{array}{ll}
\frac{d_{RL}^\nu}{dt^\nu}\hat{p}_0^\nu(t)=\lambda \hat{p}_0^\nu(t)\\
\frac{d_{RL}^\nu}{dt^\nu}\hat{p}_k^\nu(t)=\lambda
\hat{p}_k^\nu(t)-\lambda\sum_{i=1}^kq_i\hat{p}_{k-i}^\nu(t)\
\mbox{for all integer}\ k\geq 1,
\end{array}\right.
\end{equation}
with the initial conditions $\hat{p}_0^\nu(0)=1$ and
$\hat{p}_k^\nu(0)=0$ for all integer $k\geq 1$.\\
\textbf{Remarks.} For $\nu=1$, both formulas
\eqref{eq:system1-nuleq1} and \eqref{eq:system1-nugeq1} reduce to
the well-known Kolmogorov forward equations (for instance one can
specialize eq. (2.3) in Chapter 14 in \cite{KarlinTaylor}). For
$k=0$, the equations in \eqref{eq:system1-nuleq1} and
\eqref{eq:system1-nugeq1} coincide with the well-known fractional
relaxation equations (see, for example, \cite{Beghin}). For $k>0$,
they can be seen as a discrete version of the fractional master
equation (see eq. (5.14) in \cite{MeerschaertNaneVellaisamy} for
the case where the common distribution of the jumps have
continuous density; see also \cite{HilferAnton}).
\end{proposition}
\noindent\emph{Proof.} In both cases ($\nu\in (0,1)$ and
$\nu\in(1,\infty)$) the initial conditions trivially hold.\\
\emph{Case} $\nu\in (0,1)$. We have to check the equations in
\eqref{eq:system1-nuleq1}. It is known (see e.g. Theorem 2.1 in
\cite{BeghinOrsingher2010}) that
\begin{equation}\label{eq:system1-nuleq1-counting-process}
\left\{\begin{array}{ll}
\frac{d_C^\nu}{dt^\nu}P(N_\lambda^\nu(t)=0)=-\lambda P(N_\lambda^\nu(t)=0)\\
\frac{d_C^\nu}{dt^\nu}P(N_\lambda^\nu(t)=k)=-\lambda
P(N_\lambda^\nu(t)=k)+\lambda P(N_\lambda^\nu(t)=k-1)\ \mbox{for
all integer}\ k\geq 1,
\end{array}\right.
\end{equation}
with the initial conditions
$$P(N_\lambda^\nu(0)=0)=1\ \mbox{and}\ P(N_\lambda^\nu(0)=k)=0\ \mbox{for all integer}\ k\geq 1.$$
The initial conditions for $(P(N_\lambda^\nu(0)=k))_{k\geq 0}$
meet the ones for $(p_k^\nu(t))_{k\geq 0}$ in the statement of the
proposition. For $k=0$ we have
$\frac{d_C^\nu}{dt^\nu}p_0^\nu(t)=-\lambda p_0^\nu(t)$ because
$p_0^\nu(t)=P(N_\lambda^\nu(t)=0)$. For $k\geq 1$ we have several
steps. Firstly, since
$$p_k^\nu(t)=\sum_{n=1}^hq_k^{*n}P(N_\lambda^\nu(t)=n)\ \mbox{for all integer}\ h\geq k,$$
by \eqref{eq:system1-nuleq1-counting-process} we get
\begin{align*}
\frac{d_C^\nu}{dt^\nu}p_k^\nu(t)=&\sum_{n=1}^kq_k^{*n}[-\lambda
P(N_\lambda^\nu(t)=n)+\lambda P(N_\lambda^\nu(t)=n-1)]\\
=&-\lambda\sum_{n=1}^kq_k^{*n}P(N_\lambda^\nu(t)=n)+\lambda\sum_{n=1}^k\left(\sum_{i=1}^kq_{k-i}^{*(n-1)}q_i\right)P(N_\lambda^\nu(t)=n-1)\\
=&-\lambda
p_k^\nu(t)+\lambda\sum_{i=1}^kq_i\sum_{n=1}^kq_{k-i}^{*(n-1)}P(N_\lambda^\nu(t)=n-1);
\end{align*}
moreover, since $q_h^{*0}=1_{\{h=0\}}$ and $q_0^{*h}=1_{\{h=0\}}$,
we get
\begin{align*}
\sum_{i=1}^kq_i\sum_{n=1}^kq_{k-i}^{*(n-1)}P(N_\lambda^\nu(t)=n-1)=&\sum_{i=1}^{k-1}q_i\sum_{n=1}^kq_{k-i}^{*(n-1)}P(N_\lambda^\nu(t)=n-1)\\
&+q_k\sum_{n=1}^kq_0^{*(n-1)}P(N_\lambda^\nu(t)=n-1)\\
=&\sum_{i=1}^{k-1}q_i\sum_{n=2}^kq_{k-i}^{*(n-1)}P(N_\lambda^\nu(t)=n-1)+q_kP(N_\lambda^\nu(t)=0)\\
=&\sum_{i=1}^{k-1}q_i\sum_{j=1}^{k-1}q_{k-i}^{*j}P(N_\lambda^\nu(t)=j)+q_kp_0^\nu(t)\\
=&\sum_{i=1}^{k-1}q_ip_{k-i}^\nu(t)+q_kp_0^\nu(t)=\sum_{i=1}^kq_ip_{k-i}^\nu(t),
\end{align*}
and this completes the proof.\\
\emph{Case} $\nu\in(1,\infty)$. We have to check the equations in
\eqref{eq:system1-nugeq1} and we follow the same lines of the
previous case. The main difference concerns only the initial step,
i.e. the analogue of \eqref{eq:system1-nuleq1-counting-process}
presented above. Thus we only give some details on how to prove
that we have the following equations
$$\left\{\begin{array}{ll}
\frac{d_{RL}^\nu}{dt^\nu}P(\hat{N}_\lambda^\nu(t)=0)=\lambda P(\hat{N}_\lambda^\nu(t)=0)\\
\frac{d_{RL}^\nu}{dt^\nu}P(\hat{N}_\lambda^\nu(t)=k)=\lambda
P(\hat{N}_\lambda^\nu(t)=k)-\lambda P(\hat{N}_\lambda^\nu(t)=k-1)\
\mbox{for all integer}\ k\geq 1,
\end{array}\right.$$
with the initial conditions
$$P(\hat{N}_\lambda^\nu(0)=0)=1\ \mbox{and}\ P(\hat{N}_\lambda^\nu(0)=k)=0\ \mbox{for all integer}\ k\geq 1.$$
Thus, since
$\hat{N}_\lambda^\nu(\cdot)=N_\lambda(\mathcal{A}^{1/\nu}(\cdot))$,
we have
$$P(\hat{N}_\lambda^\nu(t)=k)=\int_0^\infty\frac{(\lambda x)^k}{k!}e^{-\lambda x}f_\nu(x,t)dx\ \mbox{for all integer}\ k\geq 0,$$
where $f_\nu(x,t):=f_{\mathcal{A}^{1/\nu}(t)}(x)$ is the density
of the random variable $\mathcal{A}^{1/\nu}(t)$. Then we get (the
first equality holds by eq. (5.17)-(5.18) in \cite{BeghinMacci}
with $\gamma=\frac{1}{\nu}$)
\begin{align*}
\frac{d_{RL}^\nu}{dt^\nu}P(\hat{N}_\lambda^\nu(t)=k)=&\int_0^\infty\frac{(\lambda
x)^k}{k!}e^{-\lambda x}\frac{\partial}{\partial
x}f_\nu(x,t)dx=\left[\frac{(\lambda x)^k}{k!}e^{-\lambda
x}f_\nu(x,t)\right]_{x=0}^{x=\infty}\\
&-\int_0^\infty\left(-\lambda P(N_\lambda(x)=k)+\lambda
P(N_\lambda(x)=k-1)\right)f_\nu(x,t)dx\\
=&\lambda P(\hat{N}_\lambda^\nu(t)=k)-\lambda
P(\hat{N}_\lambda^\nu(t)=k-1),
\end{align*}
and the proof is complete by taking into account
$P(\hat{N}_\lambda^\nu(t)=-1)=0$ for the case $k=0$. $\Box$

\paragraph{Remark.}
If $\nu\in(1,\infty)$, we can also obtain the alternative equation
$$\frac{d}{dt^\nu}\hat{p}_k^\nu(t)=-\lambda^{1/\nu}
\sum_{j=0}^k(-1)^j\bincoeff{1/\nu}{j}\sum_{l=0}^\infty
q_k^{*(l+j)}P(\hat{N}_\lambda^\nu(t)=l)\ \mbox{for all integer}\
k\geq 1,$$ by taking the classical derivative with respect to $t$
in
$\hat{p}_k^\nu(t)=\sum_{n=1}^kq_k^{*n}P(\hat{N}_\lambda^\nu(t)=n)$,
and by \eqref{eq:OP}. It is easy to check that, for $\nu=1$, we
meet the classical Kolmogorov equation.\\

Now let $N_1(\cdot)$ be a non-fractional Poisson process with
intensity 1, independent of $\mathcal{L}^\nu(\cdot)$ (for
$\nu\in(0,1)$) and of $\mathcal{A}^{1/\nu}(\cdot)$ (for
$\nu\in(1,\infty)$). Moreover let $S(\cdot)$ be a subordinator,
independent of all the other processes such that the following
conditions hold.
\begin{enumerate}
\item If $\nu\in(0,1)$, $M^\nu(t)$ is distributed as
$N_1(S(\mathcal{L}^\nu(t)))$; in what follows we set
$S^\nu(\cdot):=S(\mathcal{L}^\nu(\cdot))$, and we have
\begin{equation}\label{eq:mixed-densities-nuleq1}
p_k^\nu(t)=\mathbb{E}\left[\frac{(S^\nu(t))^k}{k!}e^{-S^\nu(t)}\right]=
\mathbb{E}\left[\frac{(S(\mathcal{L}^\nu(t)))^k}{k!}e^{-S(\mathcal{L}^\nu(t))}\right]\
\mbox{for all integer}\ k\geq 0.
\end{equation}
\item If $\nu\in(1,\infty)$, $\hat{M}^\nu(t)$ is distributed as
$N_1(S(\mathcal{A}^{1/\nu}(t)))$; in what follows we set
$\hat{S}^\nu(\cdot):=S(\mathcal{A}^{1/\nu}(\cdot))$, and we have
\begin{equation}\label{eq:mixed-densities-nugeq1}
\hat{p}_k^\nu(t)=\mathbb{E}\left[\frac{(\hat{S}^\nu(t))^k}{k!}e^{-\hat{S}^\nu(t)}\right]=
\mathbb{E}\left[\frac{(S(\mathcal{A}^{1/\nu}(t)))^k}{k!}e^{-S(\mathcal{A}^{1/\nu}(t))}\right]\
\mbox{for all integer}\ k\geq 0.
\end{equation}
\end{enumerate}
In such a case we talk about of \emph{mixed representation} of the
fractional versions of the process $M(\cdot)$ in \eqref{eq:cs0}.
Furthermore, in view of what follows, it is useful to introduce
the function $\kappa_S$ defined by
$\kappa_S(\theta):=\log\mathbb{E}[e^{\theta S(1)}]$ for all
$\theta\in\mathbb{R}$ (note that $\kappa_S(\cdot)$ is a
nondecreasing function and $\kappa_S(0)=0$). This function has a
strict relationship with the probability generating function $g_q$
presented above; actually, by considering standard computations on
independent random time-changes for L\'{e}vy processes (see e.g.
\cite{Bertoin}) one can easily check that
\begin{equation}\label{eq:mgf-subordinator}
\mathbb{E}\left[e^{\theta S(t)}\right]=e^{\lambda
t(g_q(1+\theta)-1)}
\end{equation}
for all $\theta$ such that $1+\theta$ belongs to the domain of
convergence of $g_q$.

We conclude with the following further formulas: if $\nu\in(0,1)$,
by \eqref{eq:Laplace-basic} we have
\begin{equation}\label{eq:*-nuleq1}
\mathbb{E}\left[e^{-\theta
S^\nu(t)}\right]=E_{\nu,1}(\kappa_S(-\theta)t^\nu)\ \mathrm{if}\
\theta\geq 0;
\end{equation}
if $\nu\in(1,\infty)$, by \eqref{eq:mgf-stable-distribution} we
have
\begin{equation}\label{eq:*-nugeq1}
\mathbb{E}\left[e^{\theta
\hat{S}^\nu(t)}\right]=\left\{\begin{array}{ll}
\exp(-(-\kappa_S(\theta))^{1/\nu}t)&\ \mathrm{if}\ \theta\leq 0\\
\infty&\ \mathrm{if}\ \theta>0.
\end{array}\right.
\end{equation}

\paragraph{On the concept of overdispersion.}
It is well-known that a real valued random variable $Y$ is said to
be overdispersed if $\mathrm{Var}[Y]-\mathbb{E}[Y]>0$; similarly,
for a process $Y(\cdot)$, we have overdispersion if all the random
variables $\{Y(t):t>0\}$ are overdispersed. Typically the compound
Poisson process $M(\cdot)$ in \eqref{eq:cs0} exhibits
overdispersion when we exclude the trivial case where the random
jumps $\{X_n:n\geq 1\}$ are all equal to 1; actually, in such a
case, we have $M(\cdot)=N_\lambda(\cdot)$, i.e. $M(\cdot)$ is a
non-fractional Poisson process, and we have
$\mathrm{Var}[M(t)]-\mathbb{E}[M(t)]=0$ for all $t>0$.

Here we want to study the same feature for the process
$M^\nu(\cdot)$ in \textbf{(FV1)}; we do not deal with
$\hat{M}^\nu(\cdot)$ in \textbf{(FV2)} because each random
variable $\hat{M}^\nu(t)$ has infinite mean (actually the moment
generating function of $\mathcal{A}^{1/\nu}(t)$ is not finite in a
neighborhood of the origin; see formula
\eqref{eq:mgf-stable-distribution} presented above). It is known
(see e.g. formulas (2.7) and (2.8) in \cite{BeghinOrsingher2009})
that
$$\left\{\begin{array}{l}
\mathbb{E}[N_\lambda^\nu(t)]=\frac{\lambda t^\nu}{\Gamma(\nu+1)}\\
\mathrm{Var}[N_\lambda^\nu(t)]=\frac{\lambda
t^\nu}{\Gamma(\nu+1)}+(\lambda t^\nu)^2Z(\nu),\ \mathrm{where}\
Z(\nu):=\frac{1}{\nu}\left(\frac{1}{\Gamma(2\nu)}-\frac{1}{\nu\Gamma^2(\nu)}\right).
\end{array}\right.$$
Thus, since we have $Z(\nu)>0$ for all $\nu\in (0,1)$ and $Z(1)=0$
(see Figure 1), the fractional Poisson process
$N_\lambda^\nu(\cdot)$ (for $\nu\in(0,1)$) exhibits
overdispersion.

\begin{figure}[h]
\hbox to\hsize\bgroup\hss
\beginpicture
\setcoordinatesystem units <2.truein,2.truein> \setplotarea x from
0 to 1.1, y from 0 to 1.1 \axis bottom ticks withvalues {$0$}
{$1$} / short at 0 1 / /
%\axis left ticks withvalues {$L_{i+1}$} {$L_{i+1}+\delta_n^i$}
% {$L_{i+1}+\delta m_i$}
%{$U_{i+1}-\delta m_i$} {$U_{i+1}-\delta_n^i$} {$U_{i+1}$}
%/ short at .2  .4 .5 .7 .8 1
%/ /

\plot 0 0  0 1.1 /

\put {$1$} [b] at  -0.02 1

\plot   0 1 0.01                    1.01112
        0.02                    1.02138
        0.03                    1.0308
        0.04                    1.03937
        0.05                    1.0471
        0.06                    1.05399
        0.07                    1.06006
        0.08                    1.0653
        0.09                    1.06973
        0.1                     1.07336
        0.11                    1.07619
        0.12                    1.07824
        0.13                    1.07952
        0.14                    1.08004
        0.15                    1.07981
        0.16                    1.07884
        0.17                    1.07716
        0.18                    1.07477
        0.19                    1.07169
        0.2                     1.06793
        0.21                    1.06351
        0.22                    1.05845
        0.23                    1.05276
        0.24                    1.04646
        0.25                    1.03957
        0.26                    1.0321
        0.27                    1.02407
        0.28                    1.0155
        0.29                    1.00641
        0.3                     0.996814
        0.31                    0.986729
        0.32                    0.976176
        0.33                    0.965173
        0.34                    0.953737
        0.35                    0.941888
        0.36                    0.929644
        0.37                    0.917022
        0.38                    0.904042
        0.39                    0.890721
        0.4                     0.877077
        0.41                    0.863129
        0.42                    0.848894
        0.43                    0.834389
        0.44                    0.819633
        0.45                    0.804642
        0.46                    0.789434
        0.47                    0.774025
        0.48                    0.758432
        0.49                    0.742672
        0.5                     0.72676
        0.51                    0.710713
        0.52                    0.694546
        0.53                    0.678274
        0.54                    0.661912
        0.55                    0.645475
        0.56                    0.628977
        0.57                    0.612433
        0.58                    0.595855
        0.59                    0.579258
        0.6                     0.562655
        0.61                    0.546057
        0.62                    0.529479
        0.63                    0.51293
        0.64                    0.496424
        0.65                    0.479972
        0.66                    0.463584
        0.67                    0.44727
        0.68                    0.431043
        0.69                    0.41491
        0.7                     0.398882
        0.71                    0.382968
        0.72                    0.367176
        0.73                    0.351515
        0.74                    0.335994
        0.75                    0.32062
        0.76                    0.305399
        0.77                    0.290341
        0.78                    0.27545
        0.79                    0.260734
        0.8                     0.246199
        0.81                    0.23185
        0.82                    0.217692
        0.83                    0.203731
        0.84                    0.189972
        0.85                    0.176419
        0.86                    0.163076
        0.87                    0.149947
        0.88                    0.137036
        0.89                    0.124346
        0.9                     0.111879
        0.91                    0.0996396
        0.92                    0.0876288
        0.93                    0.0758493
        0.94                    0.064303
        0.95                    0.0529916
        0.96                    0.0419165
        0.97                    0.0310789
        0.98                    0.0204798
        0.99                    0.01012
        1.                      0          /

%\setdashes
%  \plot 0 0.223144  0.33333 0.0664757 /
%\plot 0.33 0   0.33  0.0680549  /

\endpicture
\hss\egroup
\vglue8pt\addtocounter{figure}{1}{\centerline{\textbf{Figure 1}:
Plot of
$Z(\nu):=\frac{1}{\nu}\left(\frac{1}{\Gamma(2\nu)}-\frac{1}{\nu\Gamma^2(\nu)}\right)$
versus $\nu$.}}
\end{figure}

\noindent Furthermore, it is known (see e.g. Appendix B.4 in
\cite{KlugmanPanjerWillmot}) that
$$\left\{\begin{array}{l}
\mathbb{E}[M^\nu(t)]=\mathbb{E}[N_\lambda^\nu(t)]\mathbb{E}[X_1]\\
\mathrm{Var}[M^\nu(t)]=\mathbb{E}[N_\lambda^\nu(t)]\mathrm{Var}[X_1]+\mathrm{Var}[N_\lambda^\nu(t)]\mathbb{E}^2[X_1];
\end{array}\right.$$
then, with some computations, we get
\begin{align*}
\mathrm{Var}[M^\nu(t)]-\mathbb{E}[M^\nu(t)]=&\mathbb{E}[N_\lambda^\nu(t)](\mathrm{Var}[X_1]-\mathbb{E}[X_1])+\mathrm{Var}[N_\lambda^\nu(t)]\mathbb{E}^2[X_1]\\
=&\frac{\lambda
t^\nu}{\Gamma(\nu+1)}(\mathbb{E}[X_1^2]-\mathbb{E}[X_1])+(\lambda
t^\nu)^2Z(\nu)\mathbb{E}^2[X_1]
\end{align*}
and $\mathbb{E}[X_1^2]-\mathbb{E}[X_1]\geq 0$ because $P(X_1\geq
1)=1$. In conclusion the compound process $M^\nu(\cdot)$ in
\textbf{(FV1)} exhibits overdispersion.

\section{Examples for the process $M(\cdot)$}\label{sec:examples}
In this section we present some examples for the process
$M(\cdot)$ in \eqref{eq:cs0} which allow to define the
corresponding fractional versions $M^\nu(\cdot)$ in \textbf{(FV1)}
and $\hat{M}^\nu(\cdot)$ in \textbf{(FV2)}. More precisely we
consider suitable choices of $(q_k)_{k\geq 1}$ and $\lambda$;
moreover, in all cases, the probability mass function
$(q_k)_{k\geq 1}$ is a \emph{zero truncated negative binomial
distribution} as in Example B.3.1.5 in Appendix B.3.1 in
\cite{KlugmanPanjerWillmot}, where $\beta$ in that reference
stands for $\frac{1-\alpha}{\alpha}$, for some $\alpha\in(0,1)$,
and $r>-1$:
$$\left.\begin{array}{llll}
\ &\ &r\neq 0&r=0\ (\mbox{the limit as}\ r\to 0)\\
\mbox{probability mass function}&\ &q_k=\frac{\bincoeff{r+k-1}{k}(1-\alpha)^k}{\alpha^{-r}-1}&q_k=-\frac{(1-\alpha)^k}{k\log\alpha}\\
\mbox{probability generating function}&\
&g_q(u)=\frac{\alpha^r}{1-\alpha^r}\frac{1-(1-u(1-\alpha))^r}{(1-u(1-\alpha))^r}&
g_q(u)=\frac{\log(1-u(1-\alpha))}{\log\alpha}
\end{array}\right.$$
(where $|u|<1/(1-\alpha)$ for the probability generating
functions). In our examples we have $r=1$ (geometric distribution)
for the Polya Aeppli process, $r=-\frac{1}{2}$ (extended truncated
negative binomial distribution) for the Poisson Inverse Gaussian
process, and $r=0$ (logarithmic distribution) for the Negative
Binomial process.

\begin{example}[Polya Aeppli process]\label{ex:PA}
We choose
$$\fbox{$\alpha=1-p$ and $r=1$ for some $p\in(0,1)$ and $\lambda>0$.}$$
The probability mass functions for the two fractional Polya Aeppli
processes can be obtained from
\eqref{eq:compound-densities-nuleq1},
\eqref{eq:compound-densities-nugeq1},
\eqref{eq:fpp-density-nu-leq-1} and
\eqref{eq:fpp-density-nu-geq-1} with
$q_k^{*n}=\bincoeff{k-1}{n-1}(1-p)^np^{k-n}$ (see e.g. Example
3.19 in \cite{KlugmanPanjerWillmot}, where $\beta$ in that
reference stands for $\frac{p}{1-p}$; another interesting
reference on this process is \cite{Minkova}).\\
Alternative formulas for the mixed representation can be obtained
from \eqref{eq:mixed-densities-nuleq1} and
\eqref{eq:mixed-densities-nugeq1} by considering $S(\cdot)$ as a
compound Poisson process with exponentially distributed summands;
more precisely, for $\mu=\lambda/p$ and $\beta=(1-p)/p$, we mean
$S(\cdot)=\sum_{k=1}^{N_\mu(\cdot)}Y_k^{(\beta)}$, where
$\{Y_k^{(\beta)}:n\geq 1\}$ are i.i.d. random variables with
density $f(x)=\beta e^{-\beta x}1_{(0,\infty)}(x)$, independent of
the (non-fractional) Poisson process $N_\mu(\cdot)$ with intensity
$\mu$ (see e.g. Subsection 11.1.2 in \cite{JohnsonKempKotz}, which
concerns a slightly more general situation with the generalized
Polya Aeppli distribution and the Tweedie distribution; in
particular we recover our case with the Tweedie distribution with
parameters $(-1,\frac{\lambda(1-p)}{p^2},\frac{1-p}{p})$).
\end{example}

Then we can prove the following result for Example \ref{ex:PA} as
a consequence of Proposition \ref{prop:Kolmogorov-equations}.

\begin{proposition}\label{prop:Kolmogorov-equations-consequence}
Let us consider the probability mass functions in
\eqref{eq:compound-densities-nuleq1} and
\eqref{eq:compound-densities-nugeq1}, where $M(\cdot)$ as in
Example \ref{ex:PA}. If $\nu\in(0,1)$, we have
\begin{equation}\label{eq:system1-nuleq1-consequence}
\left\{\begin{array}{ll}
\frac{d_C^\nu}{dt^\nu}p_0^\nu(t)=-\lambda p_0^\nu(t)\\
\frac{d_C^\nu}{dt^\nu}p_k^\nu(t)-(1-\alpha)\frac{d_C^\nu}{dt^\nu}p_{k-1}^\nu(t)=-\lambda
p_k^\nu(t)+\lambda p_{k-1}^\nu(t)\ \mbox{for all integer}\ k\geq 1
\end{array}\right.
\end{equation}
with the initial conditions $p_0^\nu(0)=1$ and $p_k^\nu(0)=0$ for
all integer $k\geq 1$.\\
If $\nu\in(1,\infty)$, we have
\begin{equation}\label{eq:system1-nugeq1-consequence}
\left\{\begin{array}{ll}
\frac{d_{RL}^\nu}{dt^\nu}\hat{p}_0^\nu(t)=\lambda \hat{p}_0^\nu(t)\\
\frac{d_{RL}^\nu}{dt^\nu}\hat{p}_k^\nu(t)-(1-\alpha)\frac{d_{RL}^\nu}{dt^\nu}\hat{p}_{k-1}^\nu(t)=\lambda
\hat{p}_k^\nu(t)-\lambda \hat{p}_{k-1}^\nu(t)\ \mbox{for all
integer}\ k\geq 1
\end{array}\right.
\end{equation}
with the initial conditions $\hat{p}_0^\nu(0)=1$ and
$\hat{p}_k^\nu(0)=0$ for all integer $k\geq 1$.\\
\textbf{Remark.} The equations in
\eqref{eq:system1-nuleq1-consequence} and
\eqref{eq:system1-nugeq1-consequence} coincide for $\nu=1$. They
are equations for the probability mass functions of the
non-fractional Polya Aeppli process, and we do not know any
references where they appear.
\end{proposition}
\noindent\emph{Proof.} Firstly, for $k=0$, we have the same
equations in Proposition \ref{prop:Kolmogorov-equations}; thus we
restrict the attention on the case $k\geq 1$. If $\nu\in(0,1)$, by
\eqref{eq:system1-nuleq1} and by taking into account that
$q_h=(1-\alpha)^{h-1}\alpha$ (for all integer $h\geq 1$), we have
\begin{align*}
\frac{d_C^\nu}{dt^\nu}p_k^\nu(t)-(1-\alpha)\frac{d_C^\nu}{dt^\nu}p_{k-1}^\nu(t)
=&-\lambda p_k^\nu(t)+\lambda\sum_{i=1}^k(1-\alpha)^{i-1}\alpha p_{k-i}^\nu(t)\\
&-(1-\alpha)\left(-\lambda p_{k-1}^\nu(t)+\lambda\sum_{i=1}^{k-1}(1-\alpha)^{i-1}\alpha p_{k-1-i}^\nu(t)\right)\\
=&-\lambda p_k^\nu(t)+\lambda\alpha p_{k-1}^\nu(t)+\lambda\sum_{i=2}^k(1-\alpha)^{i-1}\alpha p_{k-i}^\nu(t)\\
&+\lambda(1-\alpha)p_{k-1}^\nu(t)-\lambda\sum_{i=1}^{k-1}(1-\alpha)^i\alpha p_{k-1-i}^\nu(t)\\
=&-\lambda p_k^\nu(t)+\lambda p_{k-1}^\nu(t).
\end{align*}
If $\nu\in(1,\infty)$, the proof follows the same lines; we have
to consider \eqref{eq:system1-nugeq1} instead of
\eqref{eq:system1-nuleq1}, and there are suitable changes of
signs. $\Box$

\begin{example}[Poisson Inverse Gaussian process]\label{ex:PIG}
We have
$$\fbox{$\alpha=1-\frac{2\beta}{1+2\beta}$, $r=-\frac{1}{2}$ and $\lambda=\lambda_{\beta,\mu}:=\frac{\mu}{\beta}((1+2\beta)^{1/2}-1)$ for some $\beta,\mu>0$.}$$
We are not aware of any reference with a formula for the
convolution densities $q_k^{*n}$ in
\eqref{eq:convolution-densities}. So
\eqref{eq:compound-densities-nuleq1} and
\eqref{eq:compound-densities-nugeq1} do not give completely
explicit formulas for the probability mass functions of the two
fractional Poisson Inverse Gaussian processes.\\
The formulas \eqref{eq:mixed-densities-nuleq1} and
\eqref{eq:mixed-densities-nugeq1} for the mixed representation can
be obtained by considering $S(\cdot)$ as an Inverse Gaussian
process; more precisely $S(t)=Y_{\mu,\beta}(t)$ should have
density
\begin{equation}\label{eq:inverse-gaussian-density}
f(x)=\frac{\mu t}{(2\pi\beta x^3)^{1/2}}\exp\left(-\frac{(x-\mu
t)^2}{2\beta x}\right)1_{(0,\infty)}(x).
\end{equation}
Moreover (see e.g. eq. (3.39) in \cite{KlugmanPanjerWillmot}) we
have
\begin{equation}\label{eq:mgf-inverse-gaussian}
\mathbb{E}\left[e^{\theta
Y_{\mu,\beta}(t)}\right]=\left\{\begin{array}{ll}
\exp\left(-\frac{\mu}{\beta}t\left((1-2\beta\theta)^{1/2}-1\right)\right)&\ \mathrm{if}\ \theta<\frac{1}{2\beta}\\
\infty&\ \mathrm{if}\ \theta\geq\frac{1}{2\beta}
\end{array}\right.
\end{equation}
\end{example}

We remark that the processes $S^\nu(\cdot)$ (if $\nu\in(0,1)$) and
$\hat{S}^\nu(\cdot)$ (if $\nu\in(1,\infty)$) for Example
\ref{ex:PIG} can be seen as a fractional Inverse Gaussian process.
A different fractional version of this process has been defined in
\cite{KumarMeerschaertVellaisamy} by subordinating a fractional
Brownian motion to an Inverse Gaussian process (in analogy with
the so-called fractional Laplace motion; see
\cite{KozubowskiMeerschaertPodgorski}). Some properties of the
process defined in \cite{KumarMeerschaertVellaisamy} are
illustrated in \cite{KumarVellaisamy}.

\begin{example}[Negative Binomial process]\label{ex:NB}
We choose
$$\fbox{$\alpha=p$, $r=0$ and $\lambda=-\log p$ for some $p\in(0,1)$.}$$
The probability mass functions for the two fractional Negative
Binomial processes can be obtained from
\eqref{eq:compound-densities-nuleq1},
\eqref{eq:compound-densities-nugeq1},
\eqref{eq:fpp-density-nu-leq-1} and
\eqref{eq:fpp-density-nu-geq-1} with $q_k^{*n}=\frac{n!}{(-\log
p)^n}\frac{(1-p)^k|s(k,n)|}{k!}$, where
$\{|s(k,n)|:k\in\{0,1,\ldots,n\}\}$ are the \emph{unsigned
Stirling numbers of the first kind} (see e.g. Theorem 6 in
\cite{PatilWani}).\\
Alternative formulas for the mixed representation can be obtained
from \eqref{eq:mixed-densities-nuleq1} and
\eqref{eq:mixed-densities-nugeq1} by considering $S(\cdot)$ as a
Gamma process; more precisely $S(t)$ should have density
$$f(x)=\frac{(\frac{p}{1-p})^t}{\Gamma(t)}
x^{t-1}e^{-\frac{p}{1-p}x}1_{(0,\infty)}(x).$$
\end{example}

\section{Processes with two fractional parameters}\label{sec:generalized-examples}
In this section we generalize the fractional versions of the
processes in Examples \ref{ex:PA} and \ref{ex:PIG} by considering
a further fractional parameter $\eta\in(0,1)$. Throughout this
section we use the notation $N_1^{(\eta)}(\cdot)$ for the process
$N_1^{1/\eta}(\cdot)$ in \eqref{eq:fpp-density-nu-geq-1}; then,
for each example, we consider independent random time-changes of
$N_1^{(\eta)}(\cdot)$ in terms of the processes $S^\nu(\cdot)$ in
\eqref{eq:mixed-densities-nuleq1} and $\hat{S}^\nu(\cdot)$ in
\eqref{eq:mixed-densities-nugeq1}, i.e. the fractional versions of
the subordinator $S(\cdot)$ which appears in the mixed
representation of $M^\nu(\cdot)$ presented in Section
\ref{sec:representations} above. Thus we have
\begin{equation}\label{eq:pmf-nuleq1}
\mbox{$p_k^{\eta,\nu}(t):=P(M^{\eta,\nu}(t)=k)$, where
$M^{\eta,\nu}(\cdot)=N_1^{(\eta)}(S^\nu(\cdot))$, if
$\nu\in(0,1)$;}
\end{equation}
\begin{equation}\label{eq:pmf-nugeq1}
\mbox{$\hat{p}_k^{\eta,\nu}(t):=P(\hat{M}^{\eta,\nu}(t)=k)$, where
$\hat{M}^{\eta,\nu}(\cdot)=N_1^{(\eta)}(\hat{S}^\nu(\cdot))$, if
$\nu\in(1,\infty)$.}
\end{equation}
We remark that we recover the cases presented above by setting
$\eta=1$. We shall see that the governing equations for these
probability mass functions have not only a fractional time
derivative but also the fractional difference operator
$(1-B)^\eta$ as \eqref{eq:OP}.

\subsection{Generalized fractional Polya Aeppli process}\label{sub:generalizedPA}
We start with the result.

\begin{proposition}\label{prop:discrete-densities-generalized-PA}
Let $S(\cdot)=\sum_{k=1}^{N_\mu(\cdot)}Y_k^{(\beta)}$ be as in
Example \ref{ex:PA}. If $\nu\in(0,1)$, the probability mass
function in \eqref{eq:pmf-nuleq1} satisfies the equations
\begin{equation}\label{eq:system-two-parameters-nuleq1-PA}
\left\{\begin{array}{ll}
\frac{d_C^\nu}{dt^\nu}p_0^{\eta,\nu}(t)=-\frac{1}{\beta}\left(\mu+\frac{d_C^\nu}{dt^\nu}\right)p_0^{\eta,\nu}(t)\\
\frac{d_C^\nu}{dt^\nu}p_k^{\eta,\nu}(t)=-\frac{1}{\beta}\left(\mu+\frac{d_C^\nu}{dt^\nu}\right)(1-B)^\eta
p_k^{\eta,\nu}(t)\ \mbox{for all integer}\ k\geq 1,
\end{array}\right.
\end{equation}
with the initial conditions $p_0^{\eta,\nu}(0)=1$ and
$p_k^{\eta,\nu}(0)=0$ for all integer $k\geq 1$.\\
If $\nu\in(1,\infty)$, the probability mass function in
\eqref{eq:pmf-nugeq1} satisfies the equations
\begin{equation}\label{eq:system-two-parameters-nugeq1-PA}
\left\{\begin{array}{ll}
\frac{d_{RL}^\nu}{dt^\nu}\hat{p}_0^{\eta,\nu}(t)=\frac{1}{\beta}\left(\mu-\frac{d_{RL}^\nu}{dt^\nu}\right)\hat{p}_0^{\eta,\nu}(t)\\
\frac{d_{RL}^\nu}{dt^\nu}\hat{p}_k^{\eta,\nu}(t)=\frac{1}{\beta}\left(\mu-\frac{d_{RL}^\nu}{dt^\nu}\right)(1-B)^\eta
\hat{p}_k^{\eta,\nu}(t)\ \mbox{for all integer}\ k\geq 1,
\end{array}\right.
\end{equation}
with the initial conditions $\hat{p}_0^{\eta,\nu}(0)=1$ and
$\hat{p}_k^{\eta,\nu}(0)=0$ for all integer $k\geq 1$.\\
\textbf{Remark.} The equations in
\eqref{eq:system-two-parameters-nuleq1-PA} and
\eqref{eq:system-two-parameters-nugeq1-PA} coincide for $\nu=1$.
\end{proposition}
\noindent\emph{Proof.} In both cases ($\nu\in(0,1)$ and
$\nu\in(1,\infty)$) the initial conditions trivially hold.\\
\emph{Case} $\nu\in(0,1)$. We have
$$P(S^\nu(t)\in G)=P(N_\mu(\mathcal{L}^\nu(t))=0)1_G(0)+\int_G f_\nu(x,t)dx\ \mbox{(for all Borel subsets $G$ of $[0,\infty)$)},$$
where $f_\nu$ is the density of the absolutely continuous part of
the random variable $S^\nu(t)$; thus, since
$P(N_\mu(\mathcal{L}^\nu(t))=0)=E_{\nu,1}(-\mu t^\nu)$, we obtain
$$p_k^{\eta,\nu}(t)=E_{\nu,1}(-\mu t^\nu)1_{\{k=0\}}+\int_0^\infty P(N_1^{(\eta)}(\mathcal{L}^\nu(z))=k)f_\nu(z,t)dz.$$
Thus we get
\begin{align*}
\frac{d_C^\nu}{dt^\nu}p_k^{\eta,\nu}(t)=
&\left(-\mu E_{\nu,1}(-\mu t^\nu)+\frac{1}{\beta}\left(\mu+\frac{d_C^\nu}{dt^\nu}\right)f_\nu(0,t)\right)1_{\{k=0\}}\\
&+\frac{1}{\beta}\left(\mu+\frac{d_C^\nu}{dt^\nu}\right)\int_0^\infty\frac{d}{dz}\{P(N_1^{(\eta)}(\mathcal{L}^\nu(z))=k)\}f_\nu(z,t)dz
\end{align*}
by (2.4.58) in \cite{KilbasSrivastavaTrujillo}, (2.5) in
\cite{BeghinMacci} and an integration by parts (where we take into
account that $P(N_1^{(\eta)}(\mathcal{L}^\nu(0))=k)=1_{\{k=0\}}$).
The last equation reduces to
$$\frac{d_C^\nu}{dt^\nu}p_k^{\eta,\nu}(t)=
\frac{1}{\beta}\left(\mu+\frac{d_C^\nu}{dt^\nu}\right)\int_0^\infty\frac{d}{dz}\{P(N_1^{(\eta)}(\mathcal{L}^\nu(z))=k)\}f_\nu(z,t)dz$$
because
\begin{equation}\label{eq:digressione-nuleq1}
-\mu E_{\nu,1}(-\mu
t^\nu)+\frac{1}{\beta}\left(\mu+\frac{d_C^\nu}{dt^\nu}\right)f_\nu(0,t)=0.
\end{equation}
We remark that eq. \eqref{eq:digressione-nuleq1} can be checked by
inspection after noting that:
$$f_\nu(0,t)=\mu\beta t^\nu E_{\nu,\nu+1}^2(-\mu t^\nu)=\mu\beta
t^\nu\sum_{j=0}^\infty\frac{(j+1)(-\mu t^\nu)^j}{\Gamma(\nu
j+\nu+1)}=-\beta\sum_{j=0}^\infty\frac{j(-\mu t^\nu)^j}{\Gamma(\nu
j+1)}$$ by eq. (2.4) in \cite{BeghinMacci} and some computations;
$$\frac{d_C^\nu}{dt^\nu}f_\nu(0,t)=\mu\beta\sum_{j=0}^\infty\frac{(j+1)(-\mu)^j\frac{d_C^\nu}{dt^\nu}t^{\nu
j+\nu}}{\Gamma(\nu
j+\nu+1)}=\mu\beta\sum_{j=0}^\infty\frac{(j+1)(-\mu
t^\nu)^j}{\Gamma(\nu j+1)}$$ by taking into account that we have
$\frac{d_C^\nu}{dt^\nu}t^{\nu j+\nu}=\frac{\Gamma(\nu
j+\nu+1)}{\Gamma(\nu j+1)}t^{\nu j}$ by (2.1.17) in
\cite{KilbasSrivastavaTrujillo}.\\
Finally we can conclude the proof by using \eqref{eq:OP} with
$\lambda=1$, and with some computations (where we distinguish the
cases $k=0$ and $k\geq 1$).\\
\emph{Case} $\nu\in(1,\infty)$. We have
$$P(\hat{S}^\nu(t)\in G)=P(N_\mu(\mathcal{A}^{1/\nu}(t))=0)1_G(0)+\int_G \hat{f}_\nu(x,t)dx\ \mbox{(for all Borel subsets $G$ of $[0,\infty)$)},$$
where $\hat{f}_\nu$ is the density of the absolutely continuous
part of the random variable $\hat{S}^\nu(t)$; thus, since
$P(N_\mu(\mathcal{A}^{1/\nu}(t))=0)=e^{-\mu^{1/\nu}t}$, we obtain
$$\hat{p}_k^{\eta,\nu}(t)=e^{-\mu^{1/\nu}t}1_{\{k=0\}}+\int_0^\infty P(N_1^{(\eta)}(\mathcal{A}^{1/\nu}(z))=k)\hat{f}_\nu(z,t)dz.$$
We follow the same lines of the proof for the previous case and we
give only some details. Firstly we get
\begin{align*}
\frac{d_{RL}^\nu}{dt^\nu}p_k^{\eta,\nu}(t)
=&\left(\mu e^{-\mu^{1/\nu}t}-\frac{1}{\beta}\left(\mu-\frac{d_{RL}^\nu}{dt^\nu}\right)f_\nu(0,t)\right)1_{\{k=0\}}\\
&-\frac{1}{\beta}\left(\mu-\frac{d_{RL}^\nu}{dt^\nu}\right)\int_0^\infty\frac{d}{dz}\{P(N_1^{(\eta)}(\mathcal{A}^{1/\nu}(z))=k)\}\hat{f}_\nu(z,t)dz.
\end{align*}
by (2.2.15) in \cite{KilbasSrivastavaTrujillo}, (5.8) in
\cite{BeghinMacci} and an integration by parts (where we take into
account that
$P(N_1^{(\eta)}(\mathcal{A}^{1/\nu}(0))=k)=1_{\{k=0\}}$). The last
equation reduces to
$$\frac{d_{RL}^\nu}{dt^\nu}p_k^{\eta,\nu}(t)
=-\frac{1}{\beta}\left(\mu-\frac{d_{RL}^\nu}{dt^\nu}\right)\int_0^\infty\frac{d}{dz}\{P(N_1^{(\eta)}(\mathcal{A}^{1/\nu}(z))=k)\}\hat{f}_\nu(z,t)dz$$
because
\begin{equation}\label{eq:digressione-nugeq1}
\mu
e^{-\mu^{1/\nu}t}-\frac{1}{\beta}\left(\mu-\frac{d_{RL}^\nu}{dt^\nu}\right)f_\nu(0,t)=0.
\end{equation}
We remark that eq. \eqref{eq:digressione-nugeq1} can be checked by
inspection after noting that:
$f_\nu(0,t)=\frac{\beta\mu^{1/\nu}t}{\nu}e^{-\mu^{1/\nu}t}$ by eq.
(5.7) in \cite{BeghinMacci};
\begin{align*}
\frac{d_{RL}^\nu}{dt^\nu}f_\nu(0,t)=&\frac{\beta\mu^{1/\nu}}{\nu}\frac{d_{RL}^\nu}{dt^\nu}\left\{te^{-\mu^{1/\nu}t}\right\}
=\frac{\beta\mu^{1/\nu}}{\nu}\frac{d_{RL}^\nu}{dt^\nu}\left(-\frac{d}{d(\mu^{1/\nu})}e^{-\mu^{1/\nu}t}\right)\\
=&-\frac{\beta\mu^{1/\nu}}{\nu}\frac{d}{d(\mu^{1/\nu})}\frac{d_{RL}^\nu}{dt^\nu}e^{-\mu^{1/\nu}t}
=-\frac{\beta\mu^{1/\nu}}{\nu}\frac{d}{d(\mu^{1/\nu})}\left\{(\mu^{1/\nu})^\nu
e^{-\mu^{1/\nu}t}\right\}\\
=&-\frac{\beta\mu^{1/\nu}}{\nu}\left(\nu\mu\cdot\mu^{-1/\nu}e^{-\mu^{1/\nu}t}-t\mu
e^{-\mu^{1/\nu}t}\right)
\end{align*}
by (2.2.15) in \cite{KilbasSrivastavaTrujillo} (for the equality
$\frac{d_{RL}^\nu}{dt^\nu}e^{-\mu^{1/\nu}t}=(\mu^{1/\nu})^\nu
e^{-\mu^{1/\nu}t}$) and some computations.\\
Finally we conclude by using \eqref{eq:OP} with $\lambda=1$ and
with some computations (where we distinguish the cases $k=0$ and
$k\geq 1$). $\Box$\\

We already know that, for $\eta=1$, $M^{\eta,\nu}(\cdot)$ and
$\hat{M}^{\eta,\nu}(\cdot)$ coincide with the fractional versions
in \textbf{(FV1)} and \textbf{(FV2)} when we deal with the Polya
Aeppli process $M(\cdot)$ as in Example \ref{ex:PA}. Then one can
check that the results in Proposition
\ref{prop:Kolmogorov-equations-consequence} can be recovered by
considering Proposition
\ref{prop:discrete-densities-generalized-PA} with $\eta=1$.
Actually $(1-B)^1p_k^{1,\nu}(t)=p_k^{1,\nu}(t)-p_{k-1}^{1,\nu}(t)$
and, if we rearrange the terms in a suitable way, the equations
\eqref{eq:system-two-parameters-nuleq1-PA} and
\eqref{eq:system-two-parameters-nugeq1-PA} with $\eta=1$ meet
\eqref{eq:system1-nuleq1-consequence} and
\eqref{eq:system1-nugeq1-consequence}, respectively, by setting
$\mu=\frac{\lambda}{p}=\frac{\lambda}{1-\alpha}$ and
$\beta=\frac{1-p}{p}=\frac{\alpha}{1-\alpha}$ as in Example
\ref{ex:PA}.

\subsection{Generalized fractional Poisson Inverse Gaussian process}\label{sub:generalizedPIG}
From now on we consider the usual symbol $\delta(x)$ for the Dirac
delta. Furthermore we consider the process
$S(\cdot)=Y_{\mu,\beta}(\cdot)$ as in Example \ref{ex:PIG}, and we
introduce the following notation.
\begin{equation}\label{eq:density-subordinator-nuleq1}
\mbox{$f_\nu(x,t)=f_{S^\nu(t)}(x)$ is the density of the random
variable $S^\nu(t)$, for $\nu\in(0,1)$;}
\end{equation}
\begin{equation}\label{eq:density-subordinator-nugeq1}
\mbox{$\hat{f}_\nu(x,t)=f_{S^\nu(t)}(x)$ is the density of the
random variable $\hat{S}^\nu(t)$, for $\nu\in(1,\infty)$.}
\end{equation}
The main result in this subsection is Proposition
\ref{prop:discrete-densities-generalized-PIG} which provides
governing equations for the probability mass functions
$(p_k^{\eta,\nu}(t))_{k\geq 0}$ and
$(\hat{p}_k^{\eta,\nu}(t))_{k\geq 0}$ when $S(\cdot)$ is as in
Example \ref{ex:PIG}. In its proof we refer to Proposition
\ref{prop:PIG-KMVextension}, i.e. a preliminary result for $f_\nu$
and $\hat{f}_\nu$; thus, in some sense, Proposition
\ref{prop:PIG-KMVextension} concerns the fractional versions of
the Poisson Inverse Gaussian process in Example \ref{ex:PIG}, and
not its generalized version considered here.

\begin{proposition}\label{prop:PIG-KMVextension}
Let $f_\nu$ and $\hat{f}_\nu$ be the functions in
\eqref{eq:density-subordinator-nuleq1} and
\eqref{eq:density-subordinator-nugeq1}, respectively. If
$\nu\in(0,1)$, we have
\begin{equation}\label{eq:auxiliary-system-PIG-nuleq1}
\frac{\partial_C^\nu}{\partial
t^\nu}\frac{\partial_C^\nu}{\partial
t^\nu}f_\nu(x,t)-2\frac{\mu}{\beta}\frac{\partial_C^\nu}{\partial
t^\nu}f_\nu(x,t)=2\frac{\mu^2}{\beta}\frac{\partial}{\partial
x}f_\nu(x,t)
\end{equation}
with the initial conditions $f_\nu(x,0)=\delta(x)$ and
$f_\nu(0,t)=0$.\\
If $\nu\in(1,\infty)$, we have
\begin{equation}\label{eq:auxiliary-system-PIG-nugeq1}
\frac{\partial_{RL}^\nu}{\partial
t^\nu}\frac{\partial_{RL}^\nu}{\partial
t^\nu}\hat{f}_\nu(x,t)+2\frac{\mu}{\beta}\frac{\partial_{RL}^\nu}{\partial
t^\nu}\hat{f}_\nu(x,t)=2\frac{\mu^2}{\beta}\frac{\partial}{\partial
x}\hat{f}_\nu(x,t)
\end{equation}
with the initial conditions $\hat{f}_\nu(x,0)=\delta(x)$ and
$\hat{f}_\nu(0,t)=0$.\\
\textbf{Remark.} The equations in
\eqref{eq:auxiliary-system-PIG-nuleq1} and
\eqref{eq:auxiliary-system-PIG-nugeq1} coincide for $\nu=1$, and
coincide with the equation Theorem 3.1 in
\cite{KumarMeerschaertVellaisamy} (with the same initial
conditions).
\end{proposition}
\noindent\emph{Proof.} We introduce the following notation:
$\tilde{g}(\theta)$, for $\theta\geq 0$, is the Laplace transform
of $g(x)$; $g^*(\sigma)$, for $\sigma\geq 0$, is the Laplace
transform of $g(t)$ (we need the second one only for
$\nu\in(0,1)$). Furthermore we recall that, by
\eqref{eq:mgf-inverse-gaussian}, we have
$$\kappa_S(-\theta)=\log\mathbb{E}[e^{-\theta S(1)}]=-\frac{\mu}{\beta}\left((1+2\beta\theta)^{1/2}-1\right).$$
\emph{Case} $\nu\in(0,1)$. We have
$\tilde{f}_\nu(\theta,t)=E_{\nu,1}(\kappa_S(-\theta)t^\nu)$ by
\eqref{eq:*-nuleq1}, and
\begin{equation}\label{eq:KBS-2.4.58}
\frac{\partial_C^\nu}{\partial
t^\nu}\tilde{f}_\nu(\theta,t)=\kappa_S(-\theta)\tilde{f}_\nu(\theta,t)
\end{equation}
by eq. (2.4.58) in \cite{KilbasSrivastavaTrujillo}. Moreover we
have
\begin{equation}\label{eq:KBS-2.4.63}
\left(\frac{\partial_C^\nu}{\partial
t^\nu}g\right)^*(\sigma)=\sigma^\nu g^*(\sigma)-\sigma^{\nu-1}g(0)
\end{equation}
by eq. (2.4.63) in \cite{KilbasSrivastavaTrujillo}. We take the
Laplace transforms with respect to $t$ in \eqref{eq:KBS-2.4.58}
and, by taking into account \eqref{eq:KBS-2.4.63} with
$g(\cdot)=\tilde{f}_\nu(\theta,\cdot)$ for the left hand side, we
get
$$\sigma^\nu\tilde{f}_\nu^*(\theta,\sigma)-\sigma^{\nu-1}\underbrace{\tilde{f}_\nu(\theta,0)}_{=E_{\nu,1}(0)=1}=\kappa_S(-\theta)\tilde{f}_\nu^*(\theta,\sigma);$$
then we have
$$0=(\sigma^\nu-\kappa_S(-\theta))\tilde{f}_\nu^*(\theta,\sigma)-\sigma^{\nu-1}=
\left(\sigma^\nu-\frac{\mu}{\beta}+\frac{\mu}{\beta}(1+2\beta\theta)^{1/2}\right)\tilde{f}_\nu^*(\theta,\sigma)-\sigma^{\nu-1}$$
and, after some computations, we obtain
\begin{equation}\label{eq:equivalent-auxiliary-system-PIG-nuleq1}
\left(\sigma^\nu\left(\sigma^\nu-2\frac{\mu}{\beta}\right)-2\frac{\mu^2}{\beta}\theta\right)\tilde{f}_\nu^*(\theta,\sigma)
-\sigma^{\nu-1}\left(\sigma^\nu-2\frac{\mu}{\beta}+\kappa_S(-\theta)\right)=0.
\end{equation}
We complete the proof for the case $\nu\in(0,1)$ by checking that
\eqref{eq:auxiliary-system-PIG-nuleq1} yields
\eqref{eq:equivalent-auxiliary-system-PIG-nuleq1}. This will be
done in what follows. We take the Laplace transforms with respect
to $x$ in \eqref{eq:auxiliary-system-PIG-nuleq1} and, after an
integration by parts where we take into account that
$f_\nu(0,t)=0$, we have
$$\frac{\partial_C^\nu}{\partial
t^\nu}\frac{\partial_C^\nu}{\partial t^\nu}\tilde{f}_\nu(\theta,t)
-2\frac{\mu}{\beta}\frac{\partial_C^\nu}{\partial
t^\nu}\tilde{f}_\nu(\theta,t)=2\frac{\mu^2}{\beta}\theta\tilde{f}_\nu(\theta,t);$$
then we take the Laplace transforms with respect to $t$ and, by
taking into account \eqref{eq:KBS-2.4.63} with
$g(\cdot)=\frac{\partial_C^\nu}{\partial
t^\nu}\tilde{f}_\nu(\theta,\cdot)$ and with
$g(\cdot)=\tilde{f}_\nu(\theta,\cdot)$ together with some
computations, we get
$$\left(\sigma^\nu-2\frac{\mu}{\beta}\right)\left(\sigma^\nu\tilde{f}_\nu^*(\theta,\sigma)-
\sigma^{\nu-1}\underbrace{\tilde{f}_\nu(\theta,0)}_{=1}\right)
-\sigma^{\nu-1}\left.\frac{\partial_C^\nu}{\partial
t^\nu}\tilde{f}_\nu(\theta,t)\right|_{t=0}
=2\frac{\mu^2}{\beta}\theta\tilde{f}_\nu^*(\theta,\sigma);$$
finally we meet \eqref{eq:equivalent-auxiliary-system-PIG-nuleq1}
because we have $\left.\frac{\partial_C^\nu}{\partial
t^\nu}\tilde{f}_\nu(\theta,t)\right|_{t=0}=\kappa_S(-\theta)$ (by
\eqref{eq:KBS-2.4.58} with $t=0$) and $\tilde{f}_\nu(\theta,0)=1$.\\
\emph{Case} $\nu\in(1,\infty)$. We have
$\tilde{\hat{f}}_\nu(\theta,t)=\exp(-(-\kappa_S(-\theta))^{1/\nu}t)$
by \eqref{eq:*-nugeq1}, and $\frac{\partial_{RL}^\nu}{\partial
t^\nu}\tilde{\hat{f}}_\nu(\theta,t)=-\kappa_S(-\theta)\tilde{\hat{f}}_\nu(\theta,t)$
and $\frac{\partial_{RL}^\nu}{\partial
t^\nu}\frac{\partial_{RL}^\nu}{\partial
t^\nu}\tilde{\hat{f}}_\nu(\theta,t)=(\kappa_S(-\theta))^2\tilde{\hat{f}}_\nu(\theta,t)$
by eq. (2.2.15) in \cite{KilbasSrivastavaTrujillo}. Then one can
obtain the equality
\begin{equation}\label{eq:equivalent-auxiliary-system-PIG-nugeq1}
\left(\frac{\partial_{RL}^\nu}{\partial
t^\nu}\frac{\partial_{RL}^\nu}{\partial
t^\nu}\tilde{\hat{f}}_\nu\right)(\theta,t)+2\frac{\mu}{\beta}\left(\frac{\partial_{RL}^\nu}{\partial
t^\nu}\tilde{\hat{f}}_\nu\right)(\theta,t)=2\frac{\mu^2}{\beta}\theta\tilde{\hat{f}}_\nu(\theta,t)
\end{equation}
with some computations. We complete the proof showing that
\eqref{eq:auxiliary-system-PIG-nugeq1} yields
\eqref{eq:equivalent-auxiliary-system-PIG-nugeq1}. We take the
Laplace transforms with respect to $x$ in
\eqref{eq:auxiliary-system-PIG-nugeq1}, and we have
$$\left(\frac{\partial_{RL}^\nu}{\partial
t^\nu}\frac{\partial_{RL}^\nu}{\partial
t^\nu}\tilde{\hat{f}}_\nu\right)(\theta,t)+2\frac{\mu}{\beta}\left(\frac{\partial_{RL}^\nu}{\partial
t^\nu}\tilde{\hat{f}}_\nu\right)(\theta,t)=\int_0^\infty
e^{-\theta x}2\frac{\mu^2}{\beta}\frac{\partial}{\partial
x}\hat{f}_\nu(x,t)dx;$$ then we meet
\eqref{eq:equivalent-auxiliary-system-PIG-nugeq1} by considering
an integration by parts for the right hand side, where we take
into account the equality $\hat{f}_\nu(0,t)=0$. $\Box$\\

Now we are ready to prove Proposition
\ref{prop:discrete-densities-generalized-PIG}.

\begin{proposition}\label{prop:discrete-densities-generalized-PIG}
Let $(p_k^{\eta,\nu}(t))_{k\geq 0}$ and
$(\hat{p}_k^{\eta,\nu}(t))_{k\geq 0}$ be the probability mass
functions in \eqref{eq:pmf-nuleq1} and \eqref{eq:pmf-nugeq1},
respectively, with $S(\cdot)=Y_{\mu,\beta}(\cdot)$ as in Example
\ref{ex:PIG}. If $\nu\in(0,1)$, we have
\begin{equation}\label{eq:system-two-parameters-nuleq1-PIG}
\frac{d_C^\nu}{dt^\nu}\frac{d_C^\nu}{dt^\nu}p_k^{\eta,\nu}(t)-2\frac{\mu}{\beta}\frac{d_C^\nu}{dt^\nu}p_k^{\eta,\nu}(t)=
2\frac{\mu^2}{\beta}(1-B)^\eta p_k^{\eta,\nu}(t)\ \mbox{for all
integer}\ k\geq 0
\end{equation}
with the initial conditions $p_0^{\eta,\nu}(0)=1$ and
$p_k^{\eta,\nu}(0)=0$ for all integer $k\geq 1$.\\
If $\nu\in(1,\infty)$, we have
\begin{equation}\label{eq:system-two-parameters-nugeq1-PIG}
\frac{d_{RL}^\nu}{dt^\nu}\frac{d_{RL}^\nu}{dt^\nu}\hat{p}_k^{\eta,\nu}(t)+2\frac{\mu}{\beta}\frac{d_{RL}^\nu}{dt^\nu}\hat{p}_k^{\eta,\nu}(t)=
2\frac{\mu^2}{\beta}(1-B)^\eta \hat{p}_k^{\eta,\nu}(t)\ \mbox{for
all integer}\ k\geq 0
\end{equation}
with the initial conditions $\hat{p}_0^{\eta,\nu}(0)=1$ and
$\hat{p}_k^{\eta,\nu}(0)=0$ for all integer $k\geq 1$.\\
\textbf{Remark.} The equations in
\eqref{eq:system-two-parameters-nuleq1-PIG} and
\eqref{eq:system-two-parameters-nugeq1-PIG} coincide for $\nu=1$.
\end{proposition}
\noindent\emph{Proof.} In both cases ($\nu\in(0,1)$ and
$\nu\in(1,\infty)$) the initial conditions trivially hold.\\
We start with the case $\nu\in(0,1)$ and we consider the function
$f_\nu$ in \eqref{eq:density-subordinator-nuleq1}. Firstly we have
$$p_k^{\eta,\nu}(t)=\int_0^\infty P(N_1^{(\eta)}(z)=k)f_\nu(z,t)dz\ (\mbox{for all integer}\ k\geq 0).$$
Then we get
\begin{align*}
\frac{d_C^\nu}{dt^\nu}\frac{d_C^\nu}{dt^\nu}p_k^{\eta,\nu}(t)-2\frac{\mu}{\beta}\frac{d_C^\nu}{dt^\nu}p_k^{\eta,\nu}(t)
=&\int_0^\infty
P(N_1^{(\eta)}(z)=k)\left(\frac{\partial_C^\nu}{\partial
t^\nu}\frac{\partial_C^\nu}{\partial t^\nu}f_\nu(z,t)
-2\frac{\mu}{\beta}\frac{\partial_C^\nu}{\partial
t^\nu}f_\nu(z,t)\right)dz\\
=&\int_0^\infty
P(N_1^{(\eta)}(z)=k)2\frac{\mu^2}{\beta}\frac{\partial}{\partial
z}f_\nu(z,t)dz
\end{align*}
by \eqref{eq:auxiliary-system-PIG-nuleq1} in Proposition
\ref{prop:PIG-KMVextension}; then
\begin{align*}
\frac{d_C^\nu}{dt^\nu}\frac{d_C^\nu}{dt^\nu}p_k^{\eta,\nu}(t)-2\frac{\mu}{\beta}\frac{d_C^\nu}{dt^\nu}p_k^{\eta,\nu}(t)
=&-2\frac{\mu^2}{\beta}\int_0^\infty\frac{d}{dz}\{P(N_1^{(\eta)}(z)=k)\}f_\nu(z,t)dz\\
=&2\frac{\mu^2}{\beta}(1-B)^\eta\underbrace{\int_0^\infty
P(N_1^{(\eta)}(z)=k)f_\nu(z,t)dz}_{=p_k^{\eta,\nu}(t)}
\end{align*}
by an integration by parts and by taking into account
$f_\nu(0,t)=0$ (for the first equality) and by \eqref{eq:OP} with
$\lambda=1$ (for the second equality). Thus
\eqref{eq:system-two-parameters-nuleq1-PIG} holds and this
completes the proof for the case $\nu\in(0,1)$. The proof for the
case $\nu\in(1,\infty)$ follows the same lines; we get
\eqref{eq:system-two-parameters-nugeq1-PIG} by considering a
suitable change of sign and, in particular,
\eqref{eq:auxiliary-system-PIG-nugeq1} instead of
\eqref{eq:auxiliary-system-PIG-nuleq1}. $\Box$\\

We already know that, for $\eta=1$, $M^{\eta,\nu}(\cdot)$ and
$\hat{M}^{\eta,\nu}(\cdot)$ coincide with the fractional versions
in \textbf{(FV1)} and \textbf{(FV2)} when we deal with the Poisson
Inverse Gaussian process $M(\cdot)$ as in Example \ref{ex:PIG}.
Then one can check that, for $k=0$, the equations in Proposition
\ref{prop:discrete-densities-generalized-PIG} can be obtained from
the equations in Proposition \ref{prop:Kolmogorov-equations}
adapted to this process. Actually, for $\nu\in(0,1)$, we have
$$\frac{d_C^\nu}{dt^\nu}p_0^{1,\nu}(t)=-\lambda_{\beta,\mu}p_0^{1,\nu}(t)=-\frac{\mu}{\beta}((1+2\beta)^{1/2}-1)p_0^{1,\nu}(t)$$
by \eqref{eq:system1-nuleq1} in Proposition
\ref{prop:Kolmogorov-equations}; then we obtain
$$\frac{d_C^\nu}{dt^\nu}\frac{d_C^\nu}{dt^\nu}p_0^{1,\nu}(t)-2\frac{\mu}{\beta}\frac{d_C^\nu}{dt^\nu}p_0^{1,\nu}(t)=2\frac{\mu^2}{\beta}p_0^{1,\nu}(t)$$
with some computations, and this meets
\eqref{eq:system-two-parameters-nuleq1-PIG} in Proposition
\ref{prop:discrete-densities-generalized-PIG} because $(1-B)^1
p_k^{1,\nu}(t)=p_k^{1,\nu}(t)-p_{k-1}^{1,\nu}(t)$ and
$p_{-1}^{1,\nu}(t)=0$ for all $t\geq 0$. For $\nu\in(1,\infty)$ we
have similar computations with suitable changes of signs (we have
\eqref{eq:system1-nugeq1} and
\eqref{eq:system-two-parameters-nugeq1-PIG} in place of
\eqref{eq:system1-nuleq1} and
\eqref{eq:system-two-parameters-nuleq1-PIG}, respectively).

\paragraph{Acknowledgements.} We would like to thank the referee
for the careful reading of the paper and for useful suggestions
which led to an improved presentation. We also thank Barbara
Pacchiarotti for Figure 1.

\end{document}